\documentclass[11pt,reqno]{amsart}

\usepackage{amsfonts, amsthm, amsmath}
\allowdisplaybreaks[4]

\usepackage{rotating}

\usepackage{tikz}

\usepackage{graphics}

\usepackage{amssymb}

\usepackage{circledsteps}

\usepackage{amsmath}

\usepackage{amscd}

\usepackage[latin2]{inputenc}

\usepackage{t1enc}

\usepackage[mathscr]{eucal}

\usepackage{indentfirst}

\usepackage{graphicx}

\usepackage{graphics}

\usepackage{pict2e}

\usepackage{mathrsfs}

\usepackage{enumerate}

\usepackage{longtable}

\usepackage{booktabs}

\usepackage[pagebackref]{hyperref}

\hypersetup{colorlinks=true}
\usepackage{cite}
\usepackage{color}
\usepackage{epic}
\usepackage{hyperref} 
\usepackage{framed}
\usepackage{mathabx}

\numberwithin{equation}{section}
\theoremstyle{plain}

\newtheorem{theorem}{Theorem}[section]

\newtheorem{lemma}[theorem]{Lemma}

\newtheorem{corollary}[theorem]{Corollary}

\newtheorem{proposition}[theorem]{Proposition}

\theoremstyle{definition}

\newtheorem{example}[theorem]{Example}

\newtheorem{?}[theorem]{Problem}

\newtheoremstyle{named}{}{}{\itshape}{}{\bfseries}{.}{.5em}{#1\thmnote{ #3}}
\theoremstyle{named}

\newcommand{\orcid}[1]{\href{https://orcid.org/#1}{\textsuperscript{\textregistered}}}

\newcommand{\f}[1]{\ifthenelse{\equal{#1}{1}}{(q;q)_\infty}{(q^{#1};q^{#1})_{\infty}}}

\def\bN{\mathbb{N}}

\def\cS{\mathcal{S}}
\def\sP{\mathscr{P}}
\def\cD{\mathcal{D}}

\def\op{\overline{p}}
\def\og{\overline{g}}

\begin{document}
\title[Two Companion $q$-Series Identities and Their Combinatorics]{Two Companion $q$-Series Identities and Their Combinatorics}

\author[H. Li]{Haijun Li\orcid{0000-0002-8747-0237}}
\address[Haijun Li]{College of Mathematics and Statistics, Chongqing University, Chongqing 401331, P.R. China}
\email{lihaijun@cqu.edu.cn; lihaijune@163.com}

\author[C. Wang]{Chenwei Wang}
\address[Chenwei Wang]{College of Mathematics and Statistics, Chongqing University, Chongqing 401331, P.R. China}
\email{wangchenwei@stu.cqu.edu.cn}

\date{\today}

\begin{abstract}
In this paper, we evaluate two companion $q$-series recently suggested for investigation by Andrews and El Bachraoui in the concluding remarks of their work. From these evaluations, we derive the corresponding partition identities and give combinatorial proofs.
\end{abstract}

\keywords{q-series, inter partition, overpartition, companion identity, bijective combinatorics.
\newline \indent 2020 {\it Mathematics Subject Classification}. 05A17, 05A19, 11P84.}

\maketitle
\section{Introduction}\label{sec:intro}

Throughout this paper, we adopt the following standard $q$-series notation; see, for example,~\cite{GR04}.
\begin{align*}
&(a; q)_n := \prod_{k=0}^{n-1} (1 - aq^k),\quad  (a; q)_\infty := \prod_{k=0}^{\infty} (1 - aq^k), \quad |q| < 1,\\
&(a_1, \ldots, a_m; q)_n := (a_1; q)_n \cdots (a_m; q)_n, \qquad n \in \mathbb{N}_0 \cup \{\infty\},
\end{align*}
where we denote the set of positive integers as $\bN$ and define $\bN_0=\bN\cup\{0\}$. 

Recently, Andrews and El Bachraoui \cite{AEB} concluded their paper by raising two problems, the second of which has already been settled by Li \cite{Li}. We now turn our attention to the first, in which they asked for evaluations of one or both of the series
\begin{align}
\mathcal S_0(q)
&:=
\sum_{n\geq1}q^n(-q^{2n+3};q^2)_\infty(-q;q)_{n-1}^2,
\label{eq:S0-def}
\\
\mathcal S_1(q)
&:=
\sum_{n\geq1}q^{2n}(-q^{2n+3};q^2)_\infty(-q;q)_{n-1}^2.
\label{eq:S1-def}
\end{align}
Their Theorem 4(a) already evaluates the aggregate
\begin{align}\label{AB:4a}
\mathcal S_0(q)+\mathcal S_1(q)
=
\sum_{n\geq1}q^n(1+q^n)
(-q^{2n+3};q^2)_\infty(-q;q)_{n-1}^2=\frac{1}{1+q}((-q; q)^2_{\infty}-(-q; q^2)_{\infty}).
\end{align}
One of the main purposes of this paper is to obtain the following theorem to solve the problem they proposed.
\begin{theorem}\label{thm:main}
For $|q|<1$, we have
\begin{align}
\cS_0(q)&=\frac{q}{(1+q)^2}(2(-q; q)^2_{\infty}-(-q; q^2)_{\infty}),\label{AB:4a1}\\
\cS_1(q)&=\frac{1}{(1+q)^2}((1-q)(-q; q)_{\infty}^2-(-q; q^2)_{\infty}).\label{AB:4a2}
\end{align}
\end{theorem}

Thanks to identity \eqref{AB:4a}, it actually suffices to evaluate just one of \eqref{AB:4a1} and \eqref{AB:4a2}, after which the other follows immediately. Nevertheless, for the sake of a complete answer and in order to illustrate different approaches, we shall provide independent proofs of both \eqref{AB:4a1} and \eqref{AB:4a2}, without recourse to \eqref{AB:4a}.

On the other hand, we also examine these two identities from a combinatorial perspective, thereby deriving two partition identities. Before stating them, let us first recall some background on integer partitions; for further details, we refer the reader to \cite{andtp}. For a given non-negative integer $n$, a {\it partition} $\lambda$ of $n$ is a weakly decreasing list of positive integers that sum up to $n$. We write $\lambda=\lambda_{1}+\lambda_{2}+\cdots+\lambda_{m}$ with $\lambda_{1}\leq \lambda_{2}\cdots\leq\lambda_{m}$, where the {\it weight} of $\lambda$ will be denoted by $|\lambda|=n$ and each $\lambda_{i}$ is called a {\it part} of $\lambda$ for $1\leq i\leq m$. The number of parts $m$ is called the {\it length} of the partition $\lambda$ and is denoted by $\ell(\lambda)$. Denote by $\sP(n)$ the set of all partitions of $n$, while its cardinality $|\sP(n)|$ is denoted as $p(n)$ and we set $\sP:=\bigcup_{n\ge 0}\sP(n)$. An {\it overpartition} \cite{CL04} of $n$ is a partition of $n$ in which the first occurrence of each part may be overlined. Following \cite{AEB}, an overpartition is called distinct when its non-overlined parts are also distinct. An overlined and a non-overlined copy of the same integer may both occur. 

Let $\og(\pi)$ denote the greatest overlined part of an overpartition $\pi$. In the summand of $\cS_0(q)+\cS_1(1)$, the factors have the following interpretation:
\begin{itemize}
\item $q^n$ makes the overlined part $\og(\pi)=n$ mandatory;
\item the factor $(-q;q)_{n-1}$ selects the remaining overlined parts from $1,\ldots,n-1$;
\item the factor $(-q;q)_{n}$ selects non-overlined parts from $1,\ldots,n$;
\item the factor $(-q^{2n+3};q^2)_\infty$ selects distinct non-overlined odd parts greater than $2n+1$.
\end{itemize}
Hence, let $D(n)$ denote the number of overpartitions $\pi$ of $n$ into distinct parts with greatest part overlined such that each non-overlined part must be $\leq \og(\pi)$ or odd and $>2\og(\pi)+1$. Then Andrews and El Bachraoui \cite{AEB} obtained
\begin{align}
\cS_0(q)+\cS_1(q)=\sum_{n\geq 1}D(n)q^n.
\end{align}
We define
\begin{align*}
\mathcal S_0(q)=\sum_{n\geq 1}d_0(n)q^n,
\quad\text{ and }\quad
\mathcal S_1(q)=\sum_{n\geq 1}d_1(n)q^n,
\end{align*}
where $d_0(n)$ counts the overpartitions $\pi$ enumerated by $D(n)$ for which no non-overlined part equals $\og(\pi)$ and $d_1(n)$ counts those for which a non-overlined part equal to $\og(\pi)$ does occur. Moreover, we define
\begin{align*}
(-q; q)^2_{\infty}=\sum_{n\geq0}\op_d(n)q^n,\quad \text{ and }\quad (-q; q^2)_{\infty}=\sum_{n\geq 0}p_{od}(n)q^n.
\end{align*}
Here $\op_d(n)$ counts overpartitions of $n$ into distinct parts, and $p_{od}(n)$ counts partitions of $n$ into distinct odd parts. Note that we adopt the convention that $d_0(n)=d_1(n)=0$ for all $n\leq 0$ and $\op_d(n)=p_{od}(n)=0$ for all $n<0$. Then the following partition identities follow as a corollary of Theorem \ref{thm:main}.

\begin{corollary}\label{cor:pi}
For every non-negative integer $n$, we have
\begin{align}
 d_0(n)+2d_0(n-1)+d_0(n-2)
&=
2\op_{d}(n-1)-p_{od}(n-1),
\label{eq:d0}
\\
 d_1(n)+2d_1(n-1)+d_1(n-2)
&=
\op_{d}(n)-\op_{d}(n-1)-p_{od}(n).
\label{eq:d1}
\end{align}
\end{corollary}

We shall show later in the paper how these two partition identities follow as corollaries of Theorem \ref{thm:main}, and, furthermore, we shall provide combinatorial proofs of them.

The paper is organized as follows. Sections \ref{sec:11} and \ref{sec:12} present $q$-series proofs of identities \eqref{AB:4a1} and \eqref{AB:4a2}, respectively. Section \ref{sec:13} is mainly devoted to a combinatorial proof of Corollary \ref{cor:pi} along with some further results, while Section \ref{sec:conclusion} offers a brief concluding summary.

\section{Proof of \eqref{AB:4a1}}\label{sec:11}

Before completing the proof of \eqref{AB:4a1}, we need a few auxiliary identities. We begin with the well-known Euler's identity
\begin{align}\label{eq:Euler}
(-q; q)_{\infty}=\frac{1}{(q; q^2)_{\infty}}.
\end{align}
Next, Ramanujan obtained the following identity (see \cite{AB, Ramanujan})
\begin{align}
(q;q^2)_\infty(aq;q^2)_\infty
\sum_{n\geq0}
\frac{(-q;q)_n(-bq;q)_nq^n}{(aq;q^2)_{n+1}}
=(-bq;q)_\infty
\sum_{n\geq0}
\frac{(q;q^2)_n(aq;q^2)_nq^{2n}}
{(-bq;q)_{2n+1}}.
\end{align}
Here, we use its specialization obtained by taking $(a, b)\mapsto(-q^2, 1)$ and then dividing both sides by $(q; q^2)_{\infty}(-q^3; q^2)_{\infty}$, which gives
\begin{align}\label{eq:Ramanujan}
\sum_{n\geq 0}\frac{(-q; q)_{n}^2q^n}{(-q^3; q^2)_{n+1}}=\frac{(-q; q)_{\infty}}{(q; q^2)_{\infty}(-q^3; q^2)_{\infty}}\sum_{n\geq 0}\frac{(q; q^2)_{n}(-q^3; q^2)_nq^{2n}}{(-q; q)_{2n+1}}.
\end{align}

More than fifty years ago, Andrews, Subbarao and Vidyasagar \cite{ASV} proved the following identity 
\begin{align}
\sum_{n\geq 0}\frac{(a; q)_n}{(b; q)_n}q^n=\frac{q(a; q)_{\infty}}{b(b; q)_{\infty}(1-aq/b)}+\frac{1-q/b}{1-aq/b}.
\end{align}
We use its specialization obtained by setting $(a, b, q)\mapsto (q, -q^2, q^2)$, which gives
\begin{align}\label{eq:ASV}
\sum_{n\geq 0}\frac{(q; q^2)_{n}q^{2n}}{(-q^2; q^2)_{n}}=\frac{1}{1+q}\left(2-\frac{(q; q^2)_{\infty}}{(-q^2; q^2)_{\infty}}\right).
\end{align}

Now we are in a position to give a complete proof of \eqref{AB:4a1}.

\begin{proof}[Proof of \eqref{AB:4a1}]
We have
\begin{align*}
\cS_0(q)&=\sum_{n\geq 1}q^n(-q^{2n+3}; q^2)_{\infty}(-q; q)^2_{n-1}\\
&=(-q; q^2)_{\infty}\sum_{n\geq 1}\frac{q^n(-q; q)^2_{n-1}}{(-q; q^2)_{n+1}}\\
&=\frac{q(-q; q^2)_{\infty}}{1+q}\sum_{n\geq 0}\frac{(-q; q)_n^2q^n}{(-q^3; q^2)_{n+1}}\quad \text{(by $n\mapsto n+1$)}\\
&=\frac{q}{1+q}\cdot\frac{(-q; q^2)_{\infty}(-q; q)_{\infty}}{(q; q^2)_{\infty}(-q^3; q^2)_{\infty}}\sum_{n\geq 0}\frac{(q; q^2)_n(-q^3; q^2)_nq^{2n}}{(-q; q)_{2n+1}}\quad\text{(by \eqref{eq:Ramanujan})}\\
&=\frac{q(-q; q)_{\infty}}{(1+q)(q; q^2)_{\infty}}\sum_{n\geq 0}\frac{(q; q^2)_nq^{2n}}{(-q^2; q^2)_n}\quad \text{(by $(-q; q)_{2n+1}=(1+q)(-q^3; q^2)_n(-q^2; q^2)_n$)}\\
&=\frac{q(-q; q)_{\infty}}{(1+q)(q; q^2)_{\infty}}\cdot \frac{1}{1+q}\left(2-\frac{(q; q^2)_{\infty}}{(-q^2; q^2)_{\infty}}\right)\quad\text{(by \eqref{eq:ASV})}\\
&=\frac{q}{(1+q)^2}(2(-q; q)_{\infty}^2-(-q; q^2)_{\infty}).\quad\text{(by \eqref{eq:Euler})}
\end{align*}
This completes the proof.
\end{proof}

\section{Proof of \eqref{AB:4a2}}\label{sec:12}

In this section, we give an independent of \eqref{AB:4a2}. Again, we need some auxiliary identities. First, Chan and Mao \cite{CM} obtained 
\begin{align}
\sum_{n\geq0}
\frac{(x,q/x;q)_nq^n}{(z,q/z;q)_{n+1}}
=
\frac{1}{x(1-z/x)(1-q/(xz))}
+
\frac{(x,q/x;q)_\infty}
{z(1-x/z)(1-q/(xz))(z,q/z;q)_\infty}.\label{eq:CM1}
\end{align}
We need to use a special form of it. Setting $z^2=-q$, the left hand side of \eqref{eq:CM1} becomes
\begin{align*}
\sum_{n\geq 0}\frac{(x, q/x; q)_nq^n}{(z, -z; q)_{n+1}}=\sum_{n\geq 0}\frac{(x, q/x; q)_nq^n}{(-q; q^2)_{n+1}},
\end{align*}
while the right hand side becomes
\begin{align*}
\frac{x}{(x-z)(x-q/z)}+\frac{x}{(z-x)(x-q/z)}\cdot\frac{(x, q/x; q)_{\infty}}{(z, -z; q)_{\infty}}
=\frac{x}{x^2+q}\left(1-\frac{(x, q/x; q)_{\infty}}{(-q; q^2)_{\infty}}\right).
\end{align*}
Hence, we obtain 
\begin{align}\label{eq:CM}
\sum_{n\geq 0}\frac{(x, q/x; q)_nq^n}{(-q; q^2)_{n+1}}=\frac{x}{x^2+q}\left(1-\frac{(x, q/x; q)_{\infty}}{(-q; q^2)_{\infty}}\right).
\end{align}

For convenience, we set
\begin{align*}
A_n(x)=(x, q/x; q)_n,\quad \text{ and }\quad A(x)=(x, q/x; q)_{\infty}.
\end{align*}
Note that we have $A_n(-1)=2(1+q^n)(-q; q)_{n-1}^2$ and $A(-1)=2(-q; q)^2_{\infty}$. Furthermore, defining the formal differential operator
\begin{align}
\cD_x:=\frac{1}{4}+\frac{1}{2}\frac{d}{dx},
\end{align}
we then have the following proposition.

\begin{proposition}\label{pro:tool}
For every integer $n\geq 1$, we have
\begin{align}
A'_n(-1)&=(q^n-1)(-q; q)_{n-1}^2;\label{t1}\\
A'(-1)&=-(-q; q)^2_{\infty};\label{t2}\\
\left.\cD_xA_n(x)\right|_{x=-1}&=q^n(-q; q)_{n-1}^2;\label{t3}\\
\left. \cD_x\left(\frac{x}{x^2+q}\right)\right|_{x=-1}&=\frac{q-3}{4(1+q)^2}.\label{t4}
\end{align}
\end{proposition}

\begin{proof}
For \eqref{t1}, the logarithmic differentiation of 
\begin{align*}
A_n(x)=
\prod_{j=0}^{n-1}(1-xq^j)
\prod_{j=0}^{n-1}\left(1-\frac{q^{j+1}}{x}\right)
\end{align*}
gives
\begin{align*}
\frac{A_n'(x)}{A_n(x)}
&=
-\sum_{j=0}^{n-1}\frac{q^j}{1-xq^j}
+
\sum_{j=0}^{n-1}
\frac{q^{j+1}/x^2}{1-q^{j+1}/x}.
\end{align*}
At $x=-1$, this becomes
\begin{align*}
A_n'(-1)=\frac{q^n-1}{2(1+q^n)}A_n(-1)=(q^n-1)(-q; q)_{n-1}^2.
\end{align*}
For \eqref{t2}, the proof proceeds in the same way as for \eqref{t1}, so we omit the details.

For \eqref{t3}, a direct expansion using \eqref{t1} yields that
\begin{align*}
\left.\cD_xA_n(x)\right|_{x=-1}=\frac{1}{4}A_n(-1)+\frac{1}{2}A_n'(-1)=q^n(-q; q)_{n-1}^2.
\end{align*}
For \eqref{t4}, a direct expansion yields that
\begin{align*}
\left. \cD_x\left(\frac{x}{x^2+q}\right)\right|_{x=-1}=\frac{-1}{4(1+q)}+\frac{q-1}{2(1+q)^2}=\frac{q-3}{4(1+q)^2}.
\end{align*}
This completes the proof.
\end{proof}

With these preparations in place, we can now give a complete proof of \eqref{AB:4a2}.

\begin{proof}[Proof of \eqref{AB:4a2}]
For $n=0$, we have $A_0(x)=1$, so
\begin{align*}
\left. \cD_xA_0(x)\right|_{x=-1}=\frac{1}{4}.
\end{align*}
Then applying the differential operator $\cD_x$ to the left hand side of \eqref{eq:CM} and then setting $x=-1$, we obtain 
\begin{align}
\left.\cD_x\left(\sum_{n\geq 0}\frac{(x, q/x; q)_nq^n}{(-q; q^2)_{n+1}}\right)\right|_{x=-1}&=\frac{1}{1+q}\left.\cD_xA_0(x)\right|_{x=-1}+\sum_{n\geq 1}\frac{q^n}{(-q; q^2)_{n+1}}\left.\cD_xA_n(x)\right|_{x=-1}\nonumber\\
&=\frac{1}{4(1+q)}+\sum_{n\geq 0}\frac{q^{2n}(-q; q)^2_{n-1}}{(-q; q^2)_{n+1}}\quad\text{(by \eqref{t3})}\nonumber\\
&=\frac{1}{4(1+q)}+\frac{\cS_1(q)}{(-q; q^2)_{\infty}}.\label{left}
\end{align}
Similarly, applying the differential operator $\cD_x$ to the right hand side of \eqref{eq:CM} and then setting $x=-1$, we obtain
\begin{align}
&\left. \cD_x\left(\frac{x}{x^2+q}\left(1-\frac{(x, q/x; q)_{\infty}}{(-q; q^2)_{\infty}}\right)\right)\right|_{x=-1}\nonumber\\
=&\frac{q-3}{4(1+q)^2}-\left.\cD_x\left(\frac{(x, q/x; q)_{\infty}}{(-q; q^2)_{\infty}}\right)\right|_{x=-1}\quad\text{(by \eqref{t4})}\nonumber\\
=&\frac{q-3}{4(1+q)^2}+\frac{(-q; q)^2_{\infty}}{2(1+q)(-q; q^2)_{\infty}}-\left.\frac{1}{2}\frac{d}{dx}\left(\frac{x}{x^2+q}\cdot\frac{A(x)}{(-q; q^2)_{\infty}}\right)\right|_{x=-1}\nonumber\\
=&\frac{q-3}{4(1+q)^2}+\frac{1-q}{(1+q)^2}\cdot\frac{(-q; q)^2_{\infty}}{(-q; q^2)_{\infty}}.\quad \text{(by \eqref{t2})}\label{right}
\end{align}
Finally, combining \eqref{left} and \eqref{right}, we obtain
\begin{align*}
\cS_1(q)&=\frac{1}{(1+q)^2}((1-q)(-q; q)_{\infty}^2-(-q; q^2)_{\infty}).
\end{align*}
This completes the proof.
\end{proof}

\section{Combinatorial proof of Corollary \ref{cor:pi}}\label{sec:13}

In this section, we first demonstrate how Corollary \ref{cor:pi} follows from Theorem \ref{thm:main}, and then present a combinatorial proof of two equivalent identities, along with an explanation of their equivalence.

\begin{proof}[Proof of Corollary \ref{cor:pi}]
Multiplying \eqref{AB:4a1} by $(1+q)^2$ gives
\begin{align*}
(1+2q+q^2)\sum_{n\geq 0}d_0(n)q^n=2\sum_{n\geq 0}\op_d(n)q^{n+1}-\sum_{n\geq 0}p_{od}q^{n+1},
\end{align*}
since $d_0(0)=0$. Comparison of the coefficient of $q^n$ proves \eqref{eq:d0}.

Similarly, \eqref{AB:4a2} gives
\begin{align*}
(1+2q+q^2)\sum_{n\geq 0}d_1(n)q^n=(1-q)\sum_{n\geq 0}\op_d(n)q^n-\sum_{n\geq 0}p_{od}(n)q^n,
\end{align*}
since $d_1(0)=0$. Comparison of the coefficient of $q^n$ proves \eqref{eq:d1}.
\end{proof}
From \eqref{eq:d1}, we also obtain the following corollary: for every $n\geq 0$, we have
\begin{align*}
\op_d(n)-\op_d(n-1)\geq p_{od}(n),
\end{align*}
where the inequality follows from the nonnegativity of $d_1(n)+2d_1(n-1)+d_1(n-2)$.

Next, we combine to prove two partition identities concerning $d_0(n)$ and $d_1(n)$, which will serve the combinatorial proof of Corollary \ref{cor:pi}. We use script letters to denote the sets of partitions corresponding to the partition numbers denoted by the corresponding lowercase letters. 
Let $p_{od>1}(n)$ denote the number of partitions of $n$ into distinct odd parts all greater than $1$, we adopt the convention that $p_{od>1}(0)=1$ and $p_{od>1}(n)=0$ for negative $n$. Then we have
\begin{proposition}\label{thm:pod}
	For every integer $n$, 
	\begin{equation}
		p_{od}(n)=p_{od>1}(n)+p_{od>1}(n-1). \label{eq:pod}
	\end{equation}
\end{proposition}
\begin{proof}
	We only need to consider all non-negative $n$.
	
	For any $\pi \in \mathscr{P}_{od}(n)$, either $\pi$ contains no part equal to $1$, in which case it belongs to $\mathscr{P}_{od>1}(n)$, or $\pi$ contains a part equal to $1$, in which case the remaining partition belongs to $\mathscr{P}_{od>1}(n-1)$. 
	
	On the other hand, for any $\lambda \in \mathscr{P}_{od>1}(n)\cup\mathscr{P}_{od>1}(n-1)$, either $\lambda\in \mathscr{P}_{od}(n)$, or $\lambda$ with a part equal to $1$ adjoined lies in $\mathscr{P}_{od}(n)$.
\end{proof}
\begin{lemma}\label{lem:d0d1}
	For every integer $n$, 
	\begin{align}
		d_0(n)&=d_0(n-1)+2d_1(n-1)+p_{od>1}(n-1),\label{eq:d02} \\
		d_1(n)+\op_{d}(n+1)&=d_1(n+1)+2d_0(n+1)+p_{od>1}(n+1).\label{eq:d12}
	\end{align}
\end{lemma}

\begin{proof}[Combinatorial proof of \eqref{eq:d02}]
	We will show that there exists a bijection $\varphi$ as follows.
	\begin{equation*}
		\begin{split}
			\varphi: \mathscr{D}_0(n) &\longrightarrow \mathscr{D}_0(n-1)\sqcup \mathscr{D}_1(n-1) \sqcup \mathscr{D}_1(n-1) \sqcup \mathscr{P}_{od>1}(n-1) \\
			\pi &\longmapsto \lambda
		\end{split}
	\end{equation*}
	
	For any given $\pi$ in $\mathscr{D}_0(n)$, we will proceed in two steps and describe the corresponding operation of $\varphi$.
	\begin{description}
	\item[Step~1] There are two cases as follows:
	\begin{description}
		\item[Case~1-I] If $\pi$ contains both the overlined part $\og(\pi)-1$ and the non-overlined part $\og(\pi)-1$, with $\og(\pi)-1>0$, then there exists a least positive integer $m$ such that $\pi$ contains both the overlined part $\og(\pi)-m$ and the non-overlined part $\og(\pi)-m$, but does not contain both the overlined part $\og(\pi)-m-1$ and the non-overlined part $\og(\pi)-m-1$ with $\og(\pi)-m-1>0$. For each $1 \le k \le m$, we merge the non-overlined part $\og(\pi)-k$ and the overlined part $\og(\pi)-k+1$ into the non-overlined odd part $2\og(\pi)-2k+1$, obtaining $\hat{\pi}$. 
		
		Notice that neither $\pi$ nor $\hat{\pi}$ contains the non-overlined part $2\og(\pi)+1$. Since $\og(\pi)-m-1>0$, we have $\og(\pi)-m <2\og(\pi)-2k+1<2\og(\pi)+1$. Thus, every non-overlined part $2\og(\pi)-2k+1$ is either not contained in $\pi$ or is used to form a new part. 
		Therefore $\hat{\pi}$ is a distinct partition that does not contain both the overlined part $\og(\hat{\pi})-1$ and the non-overlined part $\og(\hat{\pi})-1$, and does not contain the non-overlined part $\og(\hat{\pi})$. Every part of $\hat{\pi}$ greater than $\og(\hat{\pi})$ is a non-overlined odd part at least $2\og(\hat{\pi})+1$.

		\item [Case~1-II] Otherwise, if $\pi$ does not contain both the overlined part $\og(\pi)-1$ and the non-overlined part $\og(\pi)-1$, then let $\hat{\pi}=\pi$. $\hat{\pi}$ does not contain the non-overlined part $\og(\hat{\pi})$. Every part of $\hat{\pi}$ greater than $\og(\hat{\pi})$ is a non-overlined odd part greater than $2\og(\hat{\pi})+1$.
	\end{description}
	\item [Step~2] Next, there are four cases.
	\begin{description}
		\item [Case~2-I] If $\hat{\pi}$ contains neither the overlined part $\og(\hat{\pi})-1$ nor the non-overlined part $\og(\hat{\pi})-1$, then $\lambda$ is obtained by decreasing the overlined part $\og(\hat{\pi})$ by $1$ in $\hat{\pi}$. 
		We see that $\lambda$ is distinct, every part of $\lambda$ greater than $\og(\lambda)$ is a non-overlined odd part
		at least $2\og(\hat{\pi})+1=2\og(\lambda)+3$, and $\lambda$ does not contain a non-overlined part equal to $\og(\hat{\pi})-1=\og(\lambda)$, we have $\lambda \in \mathscr{D}_0(n-1)$. 
		\item [Case~2-II] If $\hat{\pi}$ contains the overlined part $\og(\hat{\pi})-1$ and does not contain the non-overlined part $\og(\hat{\pi})-1$, then we replace the overlined part $\og(\hat{\pi})$ by the non-overlined part $\og(\hat{\pi})-1$, obtaining $\lambda$. We see that $\lambda$ is distinct, every part of $\lambda$ greater than $\og(\lambda)$ is a non-overlined odd part 
		that is at least $2\og(\hat{\pi})+1=2\og(\lambda)+3$, and $\lambda$ contains the non-overlined part $\og(\hat{\pi})-1=\og(\lambda)$, we regard $\lambda$ as a partition in the first copy of $\mathscr{D}_1(n-1)$ in $\mathscr{D}_1(n-1)\sqcup \mathscr{D}_1(n-1)$. 
		\item [Case~2-III] If $\hat{\pi}$ contains the non-overlined part $\og(\hat{\pi})-1$ and does not contain the overlined part $\og(\hat{\pi})-1$, then we replace the overlined part $\og(\hat{\pi})$ by the overlined part $\og(\hat{\pi})-1$, obtaining $\lambda$. For the same reason as in Case~2-II, we can regard $\lambda$ as a partition in the second copy of $\mathscr{D}_1(n-1)$ in the disjoint union $\mathscr{D}_1(n-1)\sqcup \mathscr{D}_1(n-1)$.
		\item [Case~2-IV] If $\og(\hat{\pi})=1$, then we obtain $\lambda=\varphi(\pi)$ by removing the overlined part $1$ of $\hat{\pi}$. Then $\lambda$ does not contain any overlined part and contains only distinct non-overlined parts, each of which is odd and greater than $1$. We regard $\lambda$ as a partition belongs to $\mathscr{P}_{od>1}(n-1)$
	\end{description}
	\end{description}
    We have verified the validity of the resulting overpartitions $\lambda$ obtained from the above operations in each case, confirming that they are well-defined. 
	Next, we show that $\varphi$ is injective. Since Case~1-I and 1-II in Step~1 map $\pi$ into two disjoint sets, and the cases in Step~2 map $\hat{\pi}$ into mutually disjoint sets, it suffices to show that for each case, our operation does not produce repeated partitions. 
	
	The operation in Case~1-I produces a consecutive sequence consisting of $m$ non-overlined parts. 
	\[
	2\og(\pi)-1,\ 2\og(\pi)-3,\dots,\ 2\og(\pi)-2k+1 ,\dots,\ 2\og(\hat{\pi})+1
	\]
	whereas the parts of $\pi$ that are originally greater than $\og(\pi)$ are all non-overlined, odd, and greater than $2\og(\pi)+1$. Then we have $\og(\pi)=m+\og(\hat{\pi})$. This, in turn, guarantees that the operation in Case~1-I does not produce repeated partitions. 
	It is clear that the operations in Case~1-II, as well as the other cases in Step~2, do not produce repeated partitions. Hence, $\varphi$ is injective. 

	On the other hand, we now describe its inverse map $\varphi^{-1}$. 
	For any $\lambda\in \mathscr{D}_0(n-1)\sqcup \mathscr{D}_1(n-1) \sqcup \mathscr{D}_1(n-1) \sqcup \mathscr{P}_{od>1}(n-1)$, we will proceed in two steps.
	\begin{description}
		\item[Step~1'] There are four cases as follows:
		\begin{description}
			\item[Case~1'-I] If $\lambda \in \mathscr{D}_0(n-1)$, then we add $1$ to the overlined part $\og(\lambda)$ to obtain $\hat{\lambda}$. We see that $\hat{\lambda}$ is distinct and does not contain the overlined part $\og(\hat{\lambda})-1$, the non-overlined part $\og(\hat{\lambda})-1$ or the non-overlined part $\og(\hat{\lambda})$. And every part of $\hat{\lambda}$ greater than $\og(\hat{\lambda})$ is a non-overlined odd part at least $2\og(\lambda)+3=2\og(\hat{\lambda})+1$.
			\item[Case~1'-II] If $\lambda$ belongs to the first copy of $\mathscr{D}_1(n-1)$ in $\mathscr{D}_1(n-1)\sqcup \mathscr{D}_1(n-1)$, then we add $1$ to the non-overlined part $\og(\lambda)$, obtaining an overlined part $\og(\lambda)+1$, thereby producing $\hat{\lambda}$. 
			We see that $\hat{\lambda}$ is distinct and contains the overlined part $\og(\hat{\lambda})-1$, but does not contain the non-overlined part $\og(\hat{\lambda})-1$ or the non-overlined part $\og(\hat{\lambda})$, and every part of $\hat{\lambda}$ greater than $\og(\hat{\lambda})$ is a non-overlined odd part at least $2\og(\lambda)+3=2\og(\hat{\lambda})+1$. 
			\item[Case~1'-III] If $\lambda$ belongs to the second copy of $\mathscr{D}_1(n-1)$ in $\mathscr{D}_1(n-1)\sqcup \mathscr{D}_1(n-1)$, then we add $1$ to the overlined part $\og(\lambda)$, obtaining an overlined part $\og(\lambda)+1$, thereby producing $\hat{\lambda}$. 
			We see that $\hat{\lambda}$ is distinct and contains the non-overlined part $\og(\hat{\lambda})-1$, but does not contain the overlined part $\og(\hat{\lambda})-1$ or the non-overlined part $\og(\hat{\lambda})$, and every part of $\hat{\lambda}$ greater than $\og(\hat{\lambda})$ is a non-overlined odd part at least $2\og(\lambda)+3=2\og(\hat{\lambda})+1$. 
			\item[Case~1'-IV] If $\lambda \in \mathscr{P}_{od>1}(n-1)$, then we add an overlined part $1$ to obtain $\hat{\lambda}$.
			We see that $\hat{\lambda}$ is distinct, does not contain the non-overlined part $1$, and every other part of $\hat{\lambda}$ is a non-overlined odd part at least $3$. 
		\end{description}
		\item [Step~2']  There are two cases for $\hat{\lambda}$, as follows:
		\begin{description}
			\item [Case~2'-I] If $\hat{\lambda}$ contains the non-overlined part $2\og(\hat{\lambda})+1$, then there exists a least positive integer $m$ such that $\hat{\lambda}$ contains the non-overlined part $2(\og(\hat{\lambda})+m)-1$ but not the non-overlined part $2(\og(\hat{\lambda})+m)+1$.
			For each $ 1\le k \le m$, we split all the non-overlined part $2(\og(\hat{\lambda})+k)-1$ into the overlined part $\og(\hat{\lambda})+k$ and the non-overlined part $\og(\hat{\lambda})+k-1$, obtaining $\pi=\varphi^{-1}(\lambda)$. We see that $\pi$ is distinct, $\og(\pi)=\og(\hat{\lambda})+m$, every non-overlined part greater than $\og(\pi)$ is greater than $2(\og(\hat{\lambda})+m)+1=2(\og(\pi))+1$. And $\pi$ does not contain 
			the non-overlined part $\og(\pi)=\og(\hat{\lambda})+m$ since $\og(\hat{\lambda})<\og(\pi)<2(\og(\hat{\lambda})+m)-1$. Actually, in this case, $\pi$ must contain both the overlined part $\og(\pi)-1$ and the non-overlined part $\og(\pi)-1$. Therefore, $\pi \in \mathscr{D}_0(n)$
			\item [Case~2'-II] If $\hat{\lambda}$ does not contain the non-overlined part $2\og(\hat{\lambda})+1$, then we have $\pi = \hat{\lambda}$. In this case, we see that $\pi \in\mathscr{D}_0(n)$, and  $\pi$ must not contain both the overlined part $\og(\pi)-1$ and the non-overlined part $\og(\pi)-1$. 
		\end{description}
	\end{description}
	The same arguments as above show that the map $\varphi^{-1}$ is also well-defined and is precisely the inverse map. Hence, we conclude that the map $\varphi$ is a bijection. 
\end{proof}
\begin{example}
As an example, here we present the table of the corresponding partitions $\pi \mapsto \lambda=\varphi(\pi)$ for $n=10$, as follows.
\begin{longtable}{@{}ll@{}}
\midrule
\multicolumn{2}{@{}p{\textwidth}@{}}{\textbf{The map $\pi \mapsto \lambda$ via Case~1-I in Step~1 and Case~2-I in Step~2 }(1 partitions):} \\
\midrule
$\overline{3}+\overline{4}+3\mapsto\overline{2}+7$ &  \\
\midrule
\multicolumn{2}{@{}p{\textwidth}@{}}{\textbf{The map $\pi \mapsto \lambda$ via Case~1-II in Step~1 and Case~2-I in Step~2} (30 partitions):} \\
\midrule
$\overline{1}+\overline{2}+\overline{4}+1+2\mapsto\overline{1}+\overline{2}+\overline{3}+1+2$ & $\overline{1}+\overline{2}+\overline{5}+2\mapsto\overline{1}+\overline{2}+\overline{4}+2$ \\
$\overline{1}+\overline{3}+\overline{5}+1\mapsto\overline{1}+\overline{3}+\overline{4}+1$ & $\overline{1}+\overline{5}+1+3\mapsto\overline{1}+\overline{4}+1+3$ \\
$\overline{2}+\overline{3}+\overline{5}\mapsto\overline{2}+\overline{3}+\overline{4}$ & $\overline{2}+\overline{5}+1+2\mapsto\overline{2}+\overline{4}+1+2$ \\
$\overline{2}+\overline{5}+3\mapsto\overline{2}+\overline{4}+3$ & $\overline{3}+\overline{5}+2\mapsto\overline{3}+\overline{4}+2$ \\
$\overline{5}+2+3\mapsto\overline{4}+2+3$ & $\overline{1}+\overline{2}+\overline{6}+1\mapsto\overline{1}+\overline{2}+\overline{5}+1$ \\
$\overline{1}+\overline{3}+\overline{6}\mapsto\overline{1}+\overline{3}+\overline{5}$ & $\overline{1}+\overline{6}+1+2\mapsto\overline{1}+\overline{5}+1+2$ \\
$\overline{1}+\overline{6}+3\mapsto\overline{1}+\overline{5}+3$ & $\overline{2}+\overline{6}+2\mapsto\overline{2}+\overline{5}+2$ \\
$\overline{3}+\overline{6}+1\mapsto\overline{3}+\overline{5}+1$ & $\overline{4}+\overline{6}\mapsto\overline{4}+\overline{5}$ \\
$\overline{6}+1+3\mapsto\overline{5}+1+3$ & $\overline{6}+4\mapsto\overline{5}+4$ \\
$\overline{1}+\overline{2}+\overline{7}\mapsto\overline{1}+\overline{2}+\overline{6}$ & $\overline{1}+\overline{7}+2\mapsto\overline{1}+\overline{6}+2$ \\
$\overline{2}+\overline{7}+1\mapsto\overline{2}+\overline{6}+1$ & $\overline{3}+\overline{7}\mapsto\overline{3}+\overline{6}$ \\
$\overline{7}+1+2\mapsto\overline{6}+1+2$ & $\overline{7}+3\mapsto\overline{6}+3$ \\
$\overline{1}+\overline{8}+1\mapsto\overline{1}+\overline{7}+1$ & $\overline{2}+\overline{8}\mapsto\overline{2}+\overline{7}$ \\
$\overline{8}+2\mapsto\overline{7}+2$ & $\overline{1}+\overline{9}\mapsto\overline{1}+\overline{8}$ \\
$\overline{9}+1\mapsto\overline{8}+1$ & $\overline{10}\mapsto\overline{9}$ \\
\midrule
\multicolumn{2}{@{}p{\textwidth}@{}}{\textbf{The map $\pi \mapsto \lambda$ via Case~1-II in Step~1 and Case~2-II in Step~2} (7 partitions):} \\
\midrule
$\overline{1}+\overline{2}+7\mapsto\overline{1}+1+7$ & $\overline{1}+\overline{2}+\overline{3}+\overline{4}\mapsto\overline{1}+\overline{2}+\overline{3}+3$ \\
$\overline{1}+\overline{3}+\overline{4}+2\mapsto\overline{1}+\overline{3}+2+3$ & $\overline{2}+\overline{3}+\overline{4}+1\mapsto\overline{2}+\overline{3}+1+3$ \\
$\overline{3}+\overline{4}+1+2\mapsto\overline{3}+1+2+3$ & $\overline{1}+\overline{4}+\overline{5}\mapsto\overline{1}+\overline{4}+4$ \\
$\overline{4}+\overline{5}+1\mapsto\overline{4}+1+4$ &  \\
\midrule
\multicolumn{2}{@{}p{\textwidth}@{}}{\textbf{The map $\pi \mapsto \lambda$ via Case~1-II in Step~1 and Case~2-III in Step~2} (7 partitions): } \\
\midrule
$\overline{2}+1+7\mapsto\overline{1}+1+7$ & $\overline{1}+\overline{2}+\overline{4}+3\mapsto\overline{1}+\overline{2}+\overline{3}+3$ \\
$\overline{1}+\overline{4}+2+3\mapsto\overline{1}+\overline{3}+2+3$ & $\overline{2}+\overline{4}+1+3\mapsto\overline{2}+\overline{3}+1+3$ \\
$\overline{4}+1+2+3\mapsto\overline{3}+1+2+3$ & $\overline{1}+\overline{5}+4\mapsto\overline{1}+\overline{4}+4$ \\
$\overline{5}+1+4\mapsto\overline{4}+1+4$ &  \\
\midrule
\multicolumn{2}{@{}p{\textwidth}@{}}{\textbf{The map $\pi \mapsto \lambda$ via Case~1-II in Step~1 and Case~2-IV in Step~2} (1 partition): } \\
\midrule
$\overline{1}+9\mapsto9$ &  \\
\bottomrule
\end{longtable}
\medskip
\end{example}

\begin{proof}[Combinatorial proof of \eqref{eq:d12}]
	We will show that there exists bijections $\psi$ as follows.
	\begin{equation*}
		\begin{split}
			\psi: \mathscr{D}_1(n) \sqcup \overline{\mathscr{P}}_{d}(n+1)&\longrightarrow \mathscr{D}_1(n+1)\sqcup \mathscr{D}_0(n+1) \sqcup \mathscr{D}_0(n+1) \sqcup \mathscr{P}_{od>1}(n+1)  \\
			\pi &\longmapsto \lambda
		\end{split}
	\end{equation*}
	
	Let $g'(\pi)$ denote the greatest non-overlined part less than $\og(\pi)$ of $\pi$, let $g'(\pi)=0$ if $\pi$ contains no part less than $\og(\pi)$. Let $s(\lambda)$ denote the smallest non-overlined part greater than $2\og(\lambda)+1$, let $s(\lambda)=0$ if no part of $\lambda$ greater than $2\og(\lambda)+1$. For any given $\pi$ in $\mathscr{D}_1(n) \sqcup \overline{\mathscr{P}}_{d}(n+1)$, we divide it into the following two cases and describe the corresponding operation of $\psi$. 
	\begin{description}
		\item [Case~I] If $\pi \in \mathscr{D}_1(n)$, then $\pi$ does not contain the non-overlined part $2\og(\pi)+1$. Merge the overlined part $\og(\pi)$ and the non-overlined part $\og(\pi)$ and adding $1$ yields a non-overlined part $2\og(\pi)+1$. 
		
		Denote the resulting new partition by $\pi'$. 
		We see that $\pi'$ is distinct, all overlined parts of $\pi'$ are less than $\og(\pi)$ and all non-overlined parts of $\pi'$ that are at least $\og(\pi)$ are odd and at least $2\og(\pi)+1$. Actually, $2\og(\pi)+1$ must be contained in $\pi'$. And $\og(\pi')<\og(\pi)$, $g'(\pi)<\og(\pi)$ since $\pi \in \mathscr{D}_1(n)$.

		There is four subcases as follows:
		\begin{description}
			\item[Case~I-1] If $\og(\pi')=g'(\pi)>0$, then $\lambda=\psi(\pi)=\pi'$. In this subcase, we see that $\lambda\in \mathscr{D}_1(n+1)$ obviously. And $s(\lambda) > 0$ since the non-overlined part $2\og(\pi)+1$ which is greater than $2\og(\lambda)+1$ must be contained in $\lambda$. 
			\item[Case~I-2] If $\og(\pi')>g'(\pi)$, then $\lambda=\pi'$. In this subcase, we regard $\lambda$ as a partition in the first copy of $\mathscr{D}_0(n+1)$ in $\mathscr{D}_0(n+1) \sqcup \mathscr{D}_0(n+1)$. And  $s(\lambda)> 0$.
			\item[Case~I-3] If $\og(\pi')<g'(\pi)$, then we overline the non-overlined part $g'(\pi)$ to obtain $\lambda$. In this subcase, we regard $\lambda$ as a partition in the second copy of $\mathscr{D}_0(n+1)$ in $\mathscr{D}_0(n+1) \sqcup \mathscr{D}_0(n+1)$. And $s(\lambda)> 0$.
			\item[Case~I-4] If $\og(\pi')=g'(\pi)=0$, then $\lambda=\pi'$. At this point, $\lambda$ has no overlined part, $\lambda$ is an odd and distinct partition. All parts of $\lambda$ are greater than $1$ since $2\og(\pi)+1>1$. 
		\end{description}
		\item [Case~II] If $\pi \in \overline{\mathscr{P}}_{d}(n+1)$, let $g(\pi)$ denote the greatest non-overlined part of $\pi$, there is three subcases as follows: 
		\begin{description}
			\item[Case~II-1] If $\og(\pi)=g(\pi)$, then $\lambda=\pi$. In this case, we see that $\lambda \in \mathscr{D}_1(n+1)$ obviously. In this subcase, $s(\lambda)=0$.
			\item[Case~II-2] If $\og(\pi)>g(\pi)$, then $\lambda=\pi$. In this case, we regard $\lambda$ as a partition in the first copy of $\mathscr{D}_0(n+1)$ in $\mathscr{D}_0(n+1) \sqcup \mathscr{D}_0(n+1)$. In this subcase, $s(\lambda)=0$.
			\item[Case~II-3] If $\og(\pi)<g(\pi)$, then we overline $g(\pi)$ to obtain $\lambda$. In this case, we regard $\lambda$ as a partition in the second copy of $\mathscr{D}_0(n+1)$ in $\mathscr{D}_0(n+1) \sqcup \mathscr{D}_0(n+1)$. In this subcase, $s(\lambda)=0$.
		\end{description} 
	\end{description}
We have verified the validity of the resulting overpartitions $\lambda$ obtained from the above operations in each case, confirming that they are well-defined. 
Next, we show that $\psi$ is injective. We observe that the $\lambda$ produced in Case~1 necessarily contains a part greater than $\og(\lambda)$, whereas the $\lambda$ produced in Case~2 necessarily contains no part greater than $\og(\lambda)$. Moreover, the operations in Case~1 and Case~2 each map $\pi$ into mutually disjoint sets. Furthermore, the operation in each subcase itself evidently does not produce repeated $\lambda$'s. Hence, $\psi$ is injective. 

On the other hand, we now describe its inverse map $\psi^{-1}$. 
For any $\lambda\in \mathscr{D}_1(n+1)\sqcup \mathscr{D}_0(n+1) \sqcup \mathscr{D}_0(n+1) \sqcup \mathscr{P}_{od>1}(n+1)$, there are four cases as follows:
\begin{description}
	\item [Case~I'] If $\lambda \in \mathscr{D}_1(n+1)$, 
	then there are two subcases as follows:
	\begin{description}
		\item[Case~I'-1] If $s(\lambda)>0$, then we decrease the non-overlined part $s(\lambda)$ by $1$ and split it into two parts, the non-overlined part $(s(\lambda)-1)/2$ and the overlined part $(s(\lambda)-1)/2 $, thereby obtaining $\pi=\psi^{-1}(\lambda)$. 
		$\pi \in \mathscr{D}_1(n)$ since every part of $\pi$ greater than $\og(\pi)=(s(\lambda)-1)/2$ is odd, non-overlined and greater than $s(\lambda)=2\og(\pi)+1$. 
		\item[Case~I'-2] If $s(\lambda)=0$, then $\pi=\lambda$. 
		And $\pi\in \overline{\mathscr{P}}_{d}(n+1)$ obviously since $\mathscr{D}_1(n+1)\subset\overline{\mathscr{P}}_{d}(n+1)$.
	\end{description}
	In these two subcases, we see that $g'(\pi)=\og(\lambda)>0$.
	\item [Case~II'] If $\lambda$ belongs to the first copy of $\mathscr{D}_0(n+1)$ in $\mathscr{D}_0(n+1)\sqcup \mathscr{D}_0(n+1)$, 
	then there are two subcases as follows:
	\begin{description}
		\item[Case~II'-1] If $s(\lambda)>0$, then we decrease $s(\lambda)$ by $1$ and split it into two parts, the overlined part $(s(\lambda)-1)/2$ and the non-overlined part $(s(\lambda)-1)/2$, thereby obtaining $\pi$.
		$\pi \in \mathscr{D}_1(n)$ since every part of $\pi$ greater than $\og(\pi)=(s(\lambda)-1)/2$ is odd, non-overlined and greater than $s(\lambda)=2\og(\pi)+1$.
		\item[Case~II'-2] If $s(\lambda)=0$, then $\pi=\lambda$. 
		And $\pi\in \overline{\mathscr{P}}_{d}(n+1)$ obviously.
	\end{description}
	In these two subcases, we see that $g'(\pi)<\og(\lambda)$.
	\item [Case~III']If $\lambda$ belongs to the second copy of $\mathscr{D}_0(n+1)$ in $\mathscr{D}_0(n+1)\sqcup \mathscr{D}_0(n+1)$, 
	then we replace the overlined part $\og(\lambda)$ by the non-overlined part $\og(\lambda)$ to obtain $\lambda'$. 
	There are two subcases as follows:
	\begin{description}
		\item[Case~III'-1] If $s(\lambda)>0$, then for $\lambda'$, decrease $s(\lambda)$ by $1$ and split it into two parts, the overlined part $(s(\lambda)-1)/2$ and the non-overlined part $(s(\lambda)-1)/2$, thereby obtaining $\pi$.
		$\pi \in \mathscr{D}_1(n)$ since every part of $\pi$ greater than $\og(\pi)=(s(\lambda)-1)/2$ is odd, non-overlined and greater than $s(\lambda)=2\og(\pi)+1$.
		\item[Case~III'-2] If $s(\lambda)=0$, $\pi=\lambda'$. 
		And $\pi\in \overline{\mathscr{P}}_{d}(n+1)$ obviously.
	\end{description}
	In these two subcases, we see that $g'(\pi)>\og(\lambda')$.
	\item [Case~IV] 
	If $\lambda \in \mathscr{P}_{od>1}(n+1)$, we decrease its smallest part $\lambda_i$ by $1$ and split it into two parts, the overlined part $(\lambda_i-1)/2$ and the non-overlined part $(\lambda_i-1)/2$, thereby obtaining $\pi$. $\pi \in \mathscr{D}_1(n)$ since every part of $\pi$ greater than $\og(\pi)$ is odd, non-overlined and greater than $2\og(\pi)+1$. In this case, if we regard $\lambda$ as an overpartition without overlined part, we have $g'(\pi)=\og(\lambda)=0$.
\end{description}	
	The same arguments as above show that the map $\psi^{-1}$ is also well-defined and is precisely the inverse map. Hence, we conclude that the map $\psi$ is a bijection.
\end{proof}
\begin{example}
As an example, here we present the table of the corresponding partitions $\pi$ and $\lambda=\psi(\pi)$ for $n=8$, as follows.
\begin{longtable}{@{}ll@{}}
	\midrule
	\multicolumn{2}{@{}p{\textwidth}@{}}{\textbf{The map $\pi \mapsto \lambda$ via Case~I-1} (1 partition): $\overline{1}+\overline{3}+1+3\mapsto\overline{1}+1+7$} \\
		\midrule
	\multicolumn{2}{@{}p{\textwidth}@{}}{\textbf{The map $\pi \mapsto \lambda$ via Case~I-2} (1 partition): $\overline{2}+\overline{3}+3\mapsto\overline{2}+7$} \\
		\midrule
	\multicolumn{2}{@{}p{\textwidth}@{}}{\textbf{The map $\pi \mapsto \lambda$ via Case~I-3} (1 partition): $\overline{3}+2+3\mapsto\overline{2}+7$} \\
		\midrule
	\multicolumn{2}{@{}p{\textwidth}@{}}{\textbf{The map $\pi \mapsto \lambda$ via Case~I-4} (1 partition): $\overline{4}+4\mapsto9$} \\
		\midrule
	\multicolumn{2}{@{}p{\textwidth}@{}}{\textbf{The map $\pi \mapsto \lambda$ via Case~II-1} (6 partitions): } \\
	\midrule
	$\overline{1}+\overline{2}+\overline{3}+3\mapsto\overline{1}+\overline{2}+\overline{3}+3$ & $\overline{1}+\overline{3}+2+3\mapsto\overline{1}+\overline{3}+2+3$ \\
	$\overline{2}+\overline{3}+1+3\mapsto\overline{2}+\overline{3}+1+3$ & $\overline{3}+1+2+3\mapsto\overline{3}+1+2+3$ \\
	$\overline{1}+\overline{4}+4\mapsto\overline{1}+\overline{4}+4$ & $\overline{4}+1+4\mapsto\overline{4}+1+4$ \\
	\midrule
	\multicolumn{2}{@{}p{\textwidth}@{}}{\textbf{The map $\pi \mapsto \lambda$ via Case~II-2} (30 partitions): } \\
	\midrule
	$\overline{1}+\overline{2}+\overline{3}+1+2\mapsto\overline{1}+\overline{2}+\overline{3}+1+2$ & $\overline{1}+\overline{2}+\overline{4}+2\mapsto\overline{1}+\overline{2}+\overline{4}+2$ \\
	$\overline{1}+\overline{3}+\overline{4}+1\mapsto\overline{1}+\overline{3}+\overline{4}+1$ & $\overline{1}+\overline{4}+1+3\mapsto\overline{1}+\overline{4}+1+3$ \\
	$\overline{2}+\overline{3}+\overline{4}\mapsto\overline{2}+\overline{3}+\overline{4}$ & $\overline{2}+\overline{4}+1+2\mapsto\overline{2}+\overline{4}+1+2$ \\
	$\overline{2}+\overline{4}+3\mapsto\overline{2}+\overline{4}+3$ & $\overline{3}+\overline{4}+2\mapsto\overline{3}+\overline{4}+2$ \\
	$\overline{4}+2+3\mapsto\overline{4}+2+3$ & $\overline{1}+\overline{2}+\overline{5}+1\mapsto\overline{1}+\overline{2}+\overline{5}+1$ \\
	$\overline{1}+\overline{3}+\overline{5}\mapsto\overline{1}+\overline{3}+\overline{5}$ & $\overline{1}+\overline{5}+1+2\mapsto\overline{1}+\overline{5}+1+2$ \\
	$\overline{1}+\overline{5}+3\mapsto\overline{1}+\overline{5}+3$ & $\overline{2}+\overline{5}+2\mapsto\overline{2}+\overline{5}+2$ \\
	$\overline{3}+\overline{5}+1\mapsto\overline{3}+\overline{5}+1$ & $\overline{4}+\overline{5}\mapsto\overline{4}+\overline{5}$ \\
	$\overline{5}+1+3\mapsto\overline{5}+1+3$ & $\overline{5}+4\mapsto\overline{5}+4$ \\
	$\overline{1}+\overline{2}+\overline{6}\mapsto\overline{1}+\overline{2}+\overline{6}$ & $\overline{1}+\overline{6}+2\mapsto\overline{1}+\overline{6}+2$ \\
	$\overline{2}+\overline{6}+1\mapsto\overline{2}+\overline{6}+1$ & $\overline{3}+\overline{6}\mapsto\overline{3}+\overline{6}$ \\
	$\overline{6}+1+2\mapsto\overline{6}+1+2$ & $\overline{6}+3\mapsto\overline{6}+3$ \\
	$\overline{1}+\overline{7}+1\mapsto\overline{1}+\overline{7}+1$ & $\overline{2}+\overline{7}\mapsto\overline{2}+\overline{7}$ \\
	$\overline{7}+2\mapsto\overline{7}+2$ & $\overline{1}+\overline{8}\mapsto\overline{1}+\overline{8}$ \\
	$\overline{8}+1\mapsto\overline{8}+1$ & $\overline{9}\mapsto\overline{9}$ \\
	\midrule
	\multicolumn{2}{@{}p{\textwidth}@{}}{\textbf{The map $\pi \mapsto \lambda$ via Case~II-3} (30 partitions): } \\
	\midrule
	$1+2+6\mapsto\overline{6}+1+2$ & $1+3+5\mapsto\overline{5}+1+3$ \\
	$1+8\mapsto\overline{8}+1$ & $2+3+4\mapsto\overline{4}+2+3$ \\
	$2+7\mapsto\overline{7}+2$ & $3+6\mapsto\overline{6}+3$ \\
	$4+5\mapsto\overline{5}+4$ & $9\mapsto\overline{9}$ \\
	$\overline{1}+1+2+5\mapsto\overline{1}+\overline{5}+1+2$ & $\overline{1}+1+3+4\mapsto\overline{1}+\overline{4}+1+3$ \\
	$\overline{1}+1+7\mapsto\overline{1}+\overline{7}+1$ & $\overline{1}+2+6\mapsto\overline{1}+\overline{6}+2$ \\
	$\overline{1}+3+5\mapsto\overline{1}+\overline{5}+3$ & $\overline{1}+8\mapsto\overline{1}+\overline{8}$ \\
	$\overline{1}+\overline{2}+1+2+3\mapsto\overline{1}+\overline{2}+\overline{3}+1+2$ & $\overline{1}+\overline{2}+1+5\mapsto\overline{1}+\overline{2}+\overline{5}+1$ \\
	$\overline{1}+\overline{2}+2+4\mapsto\overline{1}+\overline{2}+\overline{4}+2$ & $\overline{1}+\overline{2}+6\mapsto\overline{1}+\overline{2}+\overline{6}$ \\
	$\overline{2}+1+2+4\mapsto\overline{2}+\overline{4}+1+2$ & $\overline{2}+1+6\mapsto\overline{2}+\overline{6}+1$ \\
	$\overline{2}+2+5\mapsto\overline{2}+\overline{5}+2$ & $\overline{2}+3+4\mapsto\overline{2}+\overline{4}+3$ \\
	$\overline{2}+7\mapsto\overline{2}+\overline{7}$ & $\overline{1}+\overline{3}+1+4\mapsto\overline{1}+\overline{3}+\overline{4}+1$ \\
	$\overline{1}+\overline{3}+5\mapsto\overline{1}+\overline{3}+\overline{5}$ & $\overline{2}+\overline{3}+4\mapsto\overline{2}+\overline{3}+\overline{4}$ \\
	$\overline{3}+1+5\mapsto\overline{3}+\overline{5}+1$ & $\overline{3}+2+4\mapsto\overline{3}+\overline{4}+2$ \\
	$\overline{3}+6\mapsto\overline{3}+\overline{6}$ & $\overline{4}+5\mapsto\overline{4}+\overline{5}$ \\
	\bottomrule
\end{longtable}
\medskip
\end{example}
Next, we can employ Lemma \ref{lem:d0d1} to obtain Corollary \ref{cor:pi}.
\begin{proof}[Another proof of Corollary \ref{cor:pi}]

For \eqref{eq:d0}, we perform the following three basic algebraic operations on the two identities in Lemma \ref{lem:d0d1}.

\begin{enumerate}
\item[(1)] Replacing $n$ by $n+1$ in \eqref{eq:d02}, we obtain
\begin{equation}
	d_0(n+1)=d_0(n)+2d_1(n)+p_{od>1}(n). \label{eq:11}
\end{equation}

\item[(2)] Replacing $n$ by $n+2$ in \eqref{eq:d02} and then swapping the two sides, we obtain
\begin{equation}
	d_0(n+1)+2d_1(n+1)+p_{od>1}(n+1)=d_0(n+2). \label{eq:12}
\end{equation}

\item[(3)] Multiplying both sides of \eqref{eq:d12} by $2$, we obtain 
\begin{equation}
	2d_1(n)+2\op_{d}(n+1)=2d_1(n+1)+4d_0(n+1)+2p_{od>1}(n+1).\label{eq:13}
\end{equation}
\end{enumerate}
Add the two sides of \eqref{eq:11}, \eqref{eq:12}, and \eqref{eq:13} correspondingly, and simplify by canceling like terms, we obtain
\begin{equation}
	2\op_{d}(n+1)=d_0(n+2)+2d_0(n+1)+d_0(n)+p_{od>1}(n+1)+p_{od>1}(n).
\end{equation}
By Proposition \ref{thm:pod}, this equation is equivalent to \eqref{eq:d0}.	

Similarly, for \eqref{eq:d1}, we perform the following three basic algebraic operations on the two identities in Lemma \ref{lem:d0d1}.

\begin{enumerate}
\item[(1)] Multiplying both sides of \eqref{eq:d02} by $2$, we obtain 
\begin{equation}
	2d_0(n)=2d_0(n-1)+4d_1(n-1)+2p_{od>1}(n-1).\label{eq:21}
\end{equation}

\item[(2)] Replacing $n$ by $n-1$ in \eqref{eq:d12}, we obtain
\begin{equation}
	d_1(n-1)+\op_{d}(n)=d_1(n)+2d_0(n)+p_{od>1}(n).\label{eq:22}
\end{equation}

\item[(3)] Replacing $n$ by $n-2$ in \eqref{eq:d12} and then swapping the two sides, we obtain
\begin{equation}
	d_1(n-1)+2d_0(n-1)+p_{od>1}(n-1)=d_1(n-2)+\op_{d}(n-1).\label{eq:23}
\end{equation}
\end{enumerate}
Add the two sides of \eqref{eq:21}, \eqref{eq:22}, and \eqref{eq:23} correspondingly, and simplify by canceling like terms, we obtain
\begin{equation}
	\op_{d}(n)=d_1(n)+2d_1(n-1)+d_1(n-2)+\op_{d}(n-1)+p_{od>1}(n)+p_{od>1}(n-1).
\end{equation}
By Proposition \ref{thm:pod}, this equation is equivalent to \eqref{eq:d1}.	
\end{proof}

\section{Conclusion}\label{sec:conclusion}
In this paper, we have completely answered the first problem proposed by Andrews and El Bachraoui \cite{AEB}. Theorem \ref{thm:main} evaluates the two companion $q$-series $\cS_0(q)$ and $\cS_1(q)$ in \eqref{eq:S0-def} and \eqref{eq:S1-def}, respectively. We gave two independent proofs: the proof of \eqref{AB:4a1} relies on a specialization of Ramanujan's identity and the Andrews, Subbarao and Vidyasagar identity, while the proof of \eqref{AB:4a2} is obtained by applying the formal differential operator $\cD_x$ to a specialization of the Chan and Mao identity. We remark that adding the two evaluations recovers the evaluation of the aggregate $\cS_0(q)+\cS_1(q)$ obtained earlier by Andrews and El Bachraoui \cite{AEB} in \eqref{AB:4a}.

Theorem \ref{thm:main} also has a combinatorial counterpart. Passing to partition numbers, we derived the two partition identities in Corollary \ref{cor:pi}, which relate the numbers $d_0(n)$ and $d_1(n)$ to the number $\op_d(n)$ of distinct overpartitions and the number $p_{od}(n)$ of partitions into distinct odd parts. We proved these identities combinatorially by constructing explicit bijections $\varphi$ and $\psi$ that prove the two identities in Lemma \ref{lem:d0d1}, from which Corollary \ref{cor:pi} follows by elementary algebra. As a byproduct, we also obtained the inequality $\op_d(n)-\op_d(n-1)\geq p_{od}(n)$ for every $n\geq 0$.

It would be of interest to find companions of other well-known $q$-series identities and to give bijective proofs of the partition identities arising from them.

\section{Declarations}
\subsection{Conflict of interest statement} On behalf of all authors, the corresponding author states that there is no conflict of interest.

\subsection{Data availability} No data was used for the research described in the article.

\subsection{Funding statement} This work received no specific grant from any funding agency in the public, commercial, or not-for-profit sectors.

\subsection{Ethics declaration} Not applicable.



\end{document}